\documentclass[11pt,psfig]{article}
\usepackage{amssymb}
\usepackage{amsmath,amsthm}
\usepackage{amsfonts}
\usepackage{graphicx}
\usepackage{graphics}
\numberwithin{equation}{section}

\usepackage{hyperref}

\begin{document} 

\begin{titlepage} 

	\centering 
	
	\scshape 
	
	\vspace*{\baselineskip} 
	

	\vspace{0.75\baselineskip} 

	\vspace{0.75\baselineskip} 

	\vspace{2\baselineskip} 
	
	
	{ \LARGE On Time-Varying Amplitude  HGARCH Model} 
	
	\vspace*{3\baselineskip} 
	

	\vspace{0.5\baselineskip} 
	
	{\scshape\large Ferdous   Mohammadi Basatini  and Saeid Rezakhah\footnote{email: rezakhah@aut.ac.ir  } %
	
	\vspace{0.7\baselineskip} 
	
	{ Amirkabir University of Technology}\\\ 424 Hafez Avenue, Tehran 15914, Iran.}  

\end{titlepage}



\vspace{0.5cm}
\begin{abstract}
The HGARCH model allows  long-memory impact in volatilities. A new  HGARCH model with time-varying amplitude is considered in this paper.  We show the stability of the model as well.  A score test is introduced  to check the time-varying behavior in amplitude.  Some value-at-risk tests  are applied to evaluate the forecastings.  Simulations are provided 
which provide further support to the proposed model. We have also have shown  the competative performance of our model in forecasting, by compairing it with HGARH and FIGARCH models for some period of SP500 indices. 
 \end{abstract}
%

Keyword: HGARCH, long-memory, time-varying, amplitude.

JEL: 13, 22, 58

\textit{Mathematics Subject Classification: 91B84, 91B30, 62F03}

\section{Introduction}
  Determining the volatility structure   is the main step in measuring risk in  financial time series. The GARCH models (Engle, 1982; Bollerslev,1986) are widely used for modeling volatility. Two kinds of structure  are recognized for GARCH models as  geometric and hyperbolic decaying that can be described as some kinds of  short-memory and long-memory respectively. Long-memory property is  present in the volatility of many
  financial data (Kwan et al., 2011). As a hyperbolic-memory model, HYGARCH (Davidson, 2004) is the most popular one  and has shown good performance in modeling long-memory behavior for many financial time series (Davidson, 2004; Tang and Shieh, 2006). The conditional variance of  HYGARCH model is a convex combination of the conditional variances of  GARCH (Bollerslev, 1986) and  FIGARCH (Baillie, 1996). The FIGARCH also shows hyperbolic-memory but has infinite variance. Li et al. (2015) argued that the conditional variance of the HYGARCH model has an unnecessarily complicated form. This motivated them to propose a new hyperbolic GARCH (HGARCH) model which is as simple as FIGARCH but has finite variance. 
  
  Financial time series often have time-varying volatilities which in many cases follow long memory in effect of exogenous and endogenous shocks. Thus models with  time-varying structure are more  appropriate for many financial time series. We consider  a HGARCH model with logistic time-varying amplitude to impose a more flexible behavior which  we call   TV-HGARCH. This time-varying amplitude allows the conditional variance to be more sensitive to  the last observation. So when a sudden shock influences the volatilities the TV-HGARCH  permits  the magnitude of variations in  the conditional variance  changes  and so make more dynamical behavior.  We show under some regularity conditions the  moments of the model  are bounded. Maximum likelihood estimators (MLEs) of the parameters  are derived. We develop a score test  to check  the presence of the time-varying amplitude  in the proposed TV-HGARCH structure.  The asymptotic behavior of MLEs and score test is verified by simulation.  Value-at-risk (VaR)  is a useful measure for quantifying the  risk which depends directly on the volatility.  The forecasts from various volatility models are evaluated and compared on the basis of how well they forecast VaR. Hence, we  perform some statistical hypothesis testing to compare  the VaR forecasts of competing models. We consider  \textit{S}\&\textit{P}500 indices from 17th February 2009 to 30th January  2015  to show the competitive behavior of TV-HGARCH model in compare to  HGARCH and FIGARCH. The paper organized as follows. The TV-HGARCH model and  the  moment properties are given in section 2. Maximum likelihood estimation is proposed in section 3. A score test is developed in Section 4 for checking time-varying amplitude. Section 5 reports the simulation studies. The VaR forecasting and its statistical testings  are provided in Section 6. The performance of the model for the empirical data of \textit{S}\&\textit{P}500 indices is reported in Section 7.  Conclusions are presented in the last section.

\section{The model}

Let $\{y_t\}$ follows a HGARCH($q,d,p$) model as 
$$y_t=\epsilon_t \sqrt{h_t} $$
\begin{equation}
h_t=\dfrac{\gamma}{\beta (1)}+w[1-\dfrac{\delta (B)}{\beta (B)}(1-B)^d]y_t^2\label{1},
\end{equation}
where  $\{\epsilon_t\}$ are identically and independently  (\textit{i.i.d.}) random variables with mean 0 and variance 1, $\gamma>0, 0<w<1$, $B$ is the back-shift operator, $\beta (x)=1-\sum_{i=1}^p\beta_ix^i$, $\delta (x)=1-\sum_{i=1}^q\delta_ix^i$  and $p, q$ are known positive integers; also  $(1-B)^d=1-\Sigma_{i=1}^{\infty}g_i B^i$ where $g_i=\dfrac{d\Gamma(i-d)}{\Gamma(1-d)\Gamma(i+1)}$ in which $0<d<1$. Let $\Upsilon_{t-1}$ be the information up to t-1 then $h_t$ is the  conditional variance as, $Var(y_t|\Upsilon_{t-1})=h_t$. The parameter $w$  is called the amplitude parameter that determines the magnitude of variations in the conditional variance (Kwan et al., 2012). For $w=1$ the model will reduce to the FIGARCH.   In this model the $h_t$ has fixed form by enriching the HGARCH model with  a time-varying  amplitude we provide a more  dynamical model for describing the volatilities.

\subsection{The Time-Varying HGARCH Model}
Let $\{y_t \}$ follows the TV-HGARCH($q,d,p$) model as
$$y_t=\epsilon_t \sqrt{h_t} $$
\begin{equation}
h_t=\dfrac{\gamma}{\beta (1)}+w_t[1-\dfrac{\delta (B)}{\beta (B)}(1-B)^d]y_t^2\label{2},
\end{equation}

where  $\{\epsilon_t\}, \gamma,  B, \beta (x), \delta (x) , p, q, (1-B)^{d}$, are defined as in (\ref{1}). Here $w_t$ is a logistic  time-varying function defined as
  \begin{eqnarray}
 w_t=\dfrac{exp(\eta \tilde{y_t})}{1+exp(\eta \tilde{y_t})}. \label{3}
 \end{eqnarray}
It is clear that $w_t$  bounded between 0,1.  $\eta>0$ is called the smoothness  parameter which determines the speed of  transition between high and low volatility. In financial time series several possible choices for the transition variable, $\tilde{y_t}$ are proposed (Dijk et al., 2002; McAleer, 2008).  We consider $\tilde{y_t}=y_{t-1}^2$ so the amplitude changes  with the size of the last observation and hence the magnitude of the last shock cause of the smooth changes of the conditional variance. \\
\subsection{Moment properties }
 Now we study the  moments  of the $\{y_t\}$. Let $\varphi=( \gamma, \beta_1,...,\beta_p, \delta_1,...,\delta_q,d)^\prime$, we can rewrite model (\ref{2}) into the form:
\begin{eqnarray*}
h_t=\phi_0+w_t \phi(B)y_t^2=\phi_0+w_t \sum_{i=1}^\infty \phi_i y_{t-i}^2,
\end{eqnarray*}
where  the $\phi_i$'s for $i=0,1,2,...$ are functions of $\varphi$. Denote $E|y_t^2|^m=M_m$, $E|\epsilon_t^2|^m=\mu_m$ and $S=\sum_{i=1}^\infty\phi_i$.  Note that $M_m=\mu_m E(h_t^m)$, assuming that $\phi_i\geq0$ and using the fact that $0\leq w_t\leq1$ it holds that
\begin{align}
h_t^m&=(\phi_0+w_t\sum_{i=1}^\infty\phi_i y_{t-i}^2)^m\\
&=\sum_{r=0}^m \displaystyle{m\choose r}\phi_0^r w_t^{m-r} \Big(\sum_{i=1}^\infty\phi_i y_{t-i}^2\Big)^{m-r}\\
&\leq\sum_{r=0}^m \displaystyle{m\choose r}\phi_0^r \Big(\sum_{i=1}^\infty\phi_i y_{t-i}^2\Big)^{m-r}\\
&=\sum_{r=0}^m \displaystyle{m\choose r}\phi_0^r \sum_{i_1=1}^\infty\sum_{i_2=1}^\infty\,...\,\sum_{i_{m-r}=1}^\infty \phi_{i_1}  \phi_{i_2}\,...\, \phi_{i_{m-r}} y_{t-i_1}^2  y_{t-i_2}^2\,...\, y_{t-i_{m-r}}^2.
\end{align}
  By the law of iterated expectations,
 \begin{eqnarray*}
 E(h_t^m)\leq \sum_{r=0}^m \displaystyle{m\choose r}\phi_0^r \sum_{i_1=1}^\infty\sum_{i_2=1}^\infty\,...\,\sum_{i_{m-r}=1}^\infty \phi_{i_1}  \phi_{i_2}\,...\, \phi_{i_{m-r}}E\Big( y_{t-i_1}^2  y_{t-i_2}^2\,...\, y_{t-i_{m-r}}^2\Big).
 \end{eqnarray*}
 Using Holder's inequality, it holds that  
 
 \begin{equation*}
 M_m \leq \mu_m\sum_{r=0}^m \displaystyle{m\choose r} \phi_0^r S^{m-r} M_{m-r}
 \end{equation*}  
and therefore
\begin{equation}
M_m\leq\dfrac{\mu_m\sum_{r=1}^m \displaystyle{m\choose r} \phi_0^r S^{m-r} M_{m-r}}{1-S^m\mu_m}. \label{03}
\end{equation}
The right hand side of (\ref{03}) is a recursive relation, so if the $M_1,M_2,...M_{m-1}$ are exist the condition
$S^m\mu_m<1$ is sufficient condition for the existence of the $M_m$. We find that for $m=6, 8$ the condition (\ref{03}) is the same as the conditions for ARFIMA-HYGARCH model presented by Kwan et al. (2012). 
As an example, we calculate the second-order moment for the TV-HGARCH(1,d,1) model
\begin{equation*}
y_t=\epsilon_t \sqrt{h_t}
\end{equation*} 
\begin{equation*}
h_t=\dfrac{\gamma}{1-\beta }+w_t[1-\dfrac{1-\delta B}{1-\beta B}(1-B)^d]y_t^2,\qquad w_t=\dfrac{exp(\eta \tilde{y_t})}{1+exp(\eta \tilde{y_t})}.\label{200}
\end{equation*}

After some calculations it holds that
\begin{equation}
h_t=\gamma+\beta h_{t-1}+w_t[{(\delta-\beta+g_1)y_{t-1}^2+\sum_{i=2}^\infty (g_i-\delta g_{i-1})}y_{t-i}^2].
\end{equation}
 Also if $(\delta-\beta+g_1)\geq0$ and  $(g_i-\delta g_{i-1}) \geq0$ for $i=2,3,...$ we have that
$$h_t \leq \gamma+\beta h_{t-1}+[{(\delta-\beta+g_1)y_{t-1}^2+\sum_{i=2}^\infty (g_i-\delta g_{i-1})}y_{t-i}^2], $$ and therefore
$$ M\leq\dfrac{\gamma}{1-(\delta+g_1)-\sum_{i=2}^\infty(g_i-\delta g_{i-1}) }.$$
Thus the $(\delta+g_1)+\sum_{i=2}^\infty(g_i-\delta g_{i-1}) <1$ is sufficient for the existence of the second-order moment of the $y_t$; i.e. $M<\infty$.

\section{Estimation}
 Let $\theta=(\varphi^\prime,\eta)^\prime$ denotes the parameter vector of the TV-HGARCH model defined in relations (\ref{2}) - (\ref{3}) and $h_t(\theta)$ refers to the conditional variance of the $y_t$ when the true parameters in TV-HGARCH model are replaced by the corresponding unknown parameters. Suppose the $y_1,...,y_T$ are a sample from the TV-HGARCH model. By assuming the normality on $\epsilon_t$, the conditional log likelihood function is $L(\theta)=-0.5 \sum_{t=1}^T l_t(\theta)$ where 
 \begin{align*}
 l_t(\theta)=\ln2\pi+\ln h_t(\theta)+\frac{y_t^2}{h_t(\theta)}.
 \end{align*}
 The derivatives of $L(\theta)$ with respect to the parameters are given as follows:
\begin{align*} \dfrac{\partial L(\theta)}{\partial \theta_{(i)}}=\sum_{t=1}^T\dfrac{1}{2h_t(\theta)}(\dfrac{y_t^2}{h_t(\theta)}-1)\dfrac{\partial h_t(\theta)}{\partial \theta_{(i)}}  
\end{align*}
where $\theta_{(i)}$ refers to the $i-th$ element of the $\theta$. The partial derivatives of $h_t(\theta)$ are obtained as:
\begin{center}
\begin{align*}
&\dfrac{\partial h_t(\theta)}{\partial \gamma}=1+\sum_{i=1}^p \beta_i \dfrac{\partial h_{t-i}}{\partial{\gamma}},\\
&\dfrac{\partial h_t(\theta)}{\partial \beta_k}=h_{t-k}+\sum_{i=1}^p \beta_i \dfrac{\partial h_{t-i}}{\partial{\beta_k}}-w_t y_{t-k}^2 \qquad k=1,..., p, \\
&\dfrac{\partial h_t(\theta)}{\partial \delta_j}=\sum_{i=1}^p \beta_i \dfrac{\partial h_{t-i}}{\partial{\delta_j}}
+w_t(1-B)^d y_{t-j}^2 \qquad j=1,..., q, \\
&\dfrac{\partial h_t(\theta)}{\partial d}=\sum_{i=1}^p \beta_i \dfrac{\partial h_{t-i}}{\partial{d}}-
w_t\delta(B)(1-B)^d \log (1-B)y_{t}^2, \\
&\dfrac{\partial h_t(\theta)}{\partial \eta}=\sum_{i=1}^p \beta_i \dfrac{\partial h_{t-i}}{\partial{\eta}}-\dfrac{\partial w_t}{\partial \eta}\Big(1-\dfrac{\delta(B)}{\beta(B)}(1-B)^d\Big)y_t^2,\\
&\dfrac{\partial w_t}{\partial \eta}=\dfrac{\tilde{y_t}exp(\eta \tilde{y_t})}{(1+exp(\eta \tilde{y_t}))^2}.
\end{align*}
\end{center}  
Here we need some numerical approaches such as  quasi-Newton algorithms  to find the  maximum likelihood estimator of the $\theta$ (Chong and Zak, 2001).

\section{Testing Time-Varying Amplitude}
For fitted HGARCH a score test is developed to check  the presence of the time-varying amplitude in the model. It is  very proper test because  only requires  the constrained estimator under $H_0$. The null hypothesis of testing time-varying amplitude corresponds to testing $H_0: \eta=0$ against $H_1:\eta>0$ in the TV-HGARCH model defined by relations (\ref{2}) - (\ref{3}). Under null hypothesis $w_t=\dfrac{1}{2}$. The null hypothesis implies the absence of the time-varying amplitude and we obtain standard  HGARCH model (Amado and Ter\"asvirta, 2008). 
Consider the conditional log-likelihood function $ L(\varphi,\eta)=-0.5\sum_{t=1}^Tl_t(\varphi,\eta)$ where
 \begin{equation*}
 l_t(\varphi,\eta)=\ln2\pi+\ln h_t(\varphi,\eta)+\frac{y_t^2}{h_t(\varphi,\eta)}.
 \end{equation*}
   At following the $\sim$ indicates the maximum likelihood estimator under $H_0$.\\                  
Let $\xi_T(\theta)=\dfrac{1}{\sqrt{T}}\sum_{t=1}^T\dfrac{\partial l_t(\theta)}{\partial\theta}$ is the average score test vector and $I(\theta)$ is the population information matrix. Consider $\theta_0=~(\varphi_0^\prime,0)^\prime$ as  true parameter vector under $H_0$. The  score test statistic  is defined as follows:
 \begin{equation}
 \lambda_{s}=\xi_T(\tilde{\theta})^\prime I^{-1}(\theta_0)\xi_T(\tilde{\theta})\sim\chi^2_{(1)} \label{5}.
 \end{equation}
Also, let $\xi_T(\theta)=(\xi_{1T}(\varphi^\prime),\xi_{2T}(\eta))^\prime $ where $\xi_{1T}(\varphi)=\dfrac{1}{\sqrt{T}}\sum_{t=1}^T\dfrac{\partial l_t(\varphi,\eta)}{\partial\varphi}$ and \\
 $\xi_{2T}(\eta)=\dfrac{1}{\sqrt{T}}\sum_{t=1}^T\dfrac{\partial l_t(\varphi,\eta)}{\partial\eta}$. So

\begin{align}
\xi_{T}(\tilde{\theta})=(0,\xi_{2T}(0))^\prime, \quad
 \xi_{2T}(0)=\dfrac{1}{\sqrt{T}}\sum_{t=1}^T\dfrac{\partial l_t(\tilde{\varphi},0) }{\partial \eta} \label{6}
\end{align}
 and
\begin{align*}
\dfrac{\partial l_t(\tilde{\varphi},0) }{\partial \eta}=(1-\frac{y_t^2}{h_t(\tilde{\varphi},0)})\dfrac{1}{h_t(\tilde{\varphi},0)}\dfrac{\partial h_t(\tilde{\varphi},0)}{\partial\eta}.
\end{align*}

Under normality, the population information matrix equals to negative expected value of the average Hessian matrix:
\begin{align*}
I(\theta) =E\left[ \dfrac{\partial^2 \log  f(y_t|\Upsilon_{t-1},\theta)} {\partial \theta \partial \theta^\prime} \right]=-E\left[ \dfrac{1}{T}\sum_{t=1}^T\dfrac{\partial^2 l_t(\theta)}{\partial \theta \partial \theta^\prime}\right]=E\left[ \dfrac{1}{T}\sum_{t=1}^T\dfrac{\partial l_t(\theta)}{\partial \theta}\dfrac{\partial l_t(\theta)}{\partial \theta^\prime}\right].   
\end{align*}

Note  (\ref{5}) depend on the unknown parameter value $\theta_0$ so it is useless. It is common to evaluate the $I^{-1}(\theta_0)$ at the $\tilde{\theta}$ to get a usable statistic. 
Hence 
\begin{align}
I(\tilde{\theta})=
\begin{bmatrix}
\tilde{I}_{11}& \tilde{I}_{12}\\
\tilde{I}_{21}& \tilde{I}_{22}  \label{7}
\end{bmatrix}
\end{align}
where
 \begin{align}
\tilde{I}_{11}=\tilde{\kappa}J,\qquad \tilde{I}_{12}=\tilde{I}_{21}=\tilde{\kappa}R, \qquad  \tilde{I}_{22}=\tilde{\kappa}Q, \label{8}
 \end{align}
\begin{align*}
& \tilde{\kappa}=\dfrac{1}{T}\sum_{t=1}^T\Big(\frac{y_t^2}{h_t(\tilde{\varphi},0)}-1\Big)^2\\ 
&J=\frac{1}{T}\sum_{t=1}^T\frac{1}{h_t^2(\tilde{\varphi},0)} \Big(\frac{\partial h_t(\tilde{\varphi},0)}{\partial\varphi}\Big) \Big(\frac{\partial h_t(\tilde{\varphi},0)}{\partial\varphi^\prime}\Big)\\
&R=\frac{1}{T}\sum_{t=1}^T\frac{1}{h_t^2(\tilde{\varphi},0)} \Big(\frac{\partial h_t(\tilde{\varphi},0)}{\partial\eta}\Big) \Big(\frac{\partial h_t(\tilde{\varphi},0)}{\partial\varphi}\Big)\\ 
&Q=\frac{1}{T}\sum_{t=1}^T\frac{1}{h_t^2(\tilde{\varphi},0)} \Big(\frac{\partial h_t(\tilde{\varphi},0)}{\partial\eta}\Big)^2.
 \end{align*}
Denote
\begin{equation}
S(\tilde{\varphi})= \xi_{2T}(0)=\dfrac{1}{\sqrt{T}}\sum_{t=1}^T\dfrac{\partial l_t(\tilde{\varphi},0) }{\partial \eta},\label{9}
\end{equation}

then by substituting  (\ref{6}) - (\ref{9}) in (\ref{6}) the score test statistic can be obtained as
\begin{equation}
\lambda_{s}=\dfrac{S^2(\tilde{\varphi})}{\tilde{\kappa}(Q-R^\prime J^{-1}R)}. \label{10}
\end{equation}  
Hence if $\tilde{\theta}= (\tilde{\varphi}^\prime,0)^\prime$  is asymptotically normal then  under $H_0:\eta=0$ the $\lambda_{s}$ will  asymptotically follows the chi-squared distribution with 1 degree of freedom under some regularity conditions (Li et al., 2011).

\section{\textbf{Simulation Study}}
This section conducts two simulation experiments to investigate  the  consistency of the MLEs (section 3) and the asymptotic behavior of the score test (section 4).  We consider three sample sizes, n=300, 500 and 1000  in  two experiments, and there are 1000 replications for each sample size. In each generated sequence the first 1000 observations have been discarded to avoid the initialization effects, so there are 1000+n observations generated each time. We simulate the data from a TV-HGARCH(1,d,1) model as follows:

$$y_t=\epsilon_t \sqrt{h_t} $$
\begin{equation}
h_t=\dfrac{\gamma}{1-\beta }+w_t[1-\dfrac{1-\delta B}{1-\beta B}(1-B)^d]y_t^2,\quad w_t=\dfrac{exp(\eta {y_{t-1}^2})}{1+exp(\eta {y_{t-1}^2})}. \label{0100}
\end{equation}

where $\{\epsilon_t \}$ are \textit{iid} standard normal variables. 

In the first experiment  the value  of the parameter vector is  
$\theta=(\gamma,\beta,\delta,d,\eta)^\prime=(.3,.4,.2,.7,1)^\prime.$
The MLE values in section 3  are calculated,  the biases (Bias) and the root mean squared error (RMSE) are summarized in Table 1. It is observed  that both Bias and RMSE are generally small and  decrease as the  sample size increases.

The second experiment is   conducted to evaluate the empirical sizes and powers of the score test statistic $\lambda_s$ in section 4. The value  of the parameter vector is  
$\theta=(\gamma,\beta,\delta,d,\eta)^\prime=(.3,.4,.2,.7,\eta)^\prime,$
when $\eta=0$ corresponds to the size and  $\eta>0$ correspond to the power of the test. We consider three different values $\eta$=0.4, 1 and 3 and two significance values .05 and .10. The empirical rejection rates are reported in  Table 2. It can be seen that the empirical sizes are all close to the nominal values and this closeness increases as the sample size increases also empirical powers are increasing function of the sample size and of the $\eta$.  
 
\begin{table}[t] 
\begin{small}
\caption{Estimation results of the TV-HGARCH model based on 1000 replications.}
\end{small}
\begin{small}
\begin{center}
\begin{tabular}{c c c c c c c c c c c c c }
\hline \hline
   & & & \multicolumn{2}{c}{n=300} & & \multicolumn{2}{c}{n=500} & & \multicolumn{3}{c}{n=1000} \\
 \cline{3-5} \cline{7-9} \cline{11-13}
 parameter & Real value  & & Bias&  RMSE & &  Bias& RMSE & &&  Bias& RMSE &\\
\hline
$\gamma$&0.3 &  & 0.031 & 0.162 & & 0.017& 0.130& && 0.005 & 0.030\\
$\beta$&0.4 &  & 0.030 & 0.014 & & 0.008 & 0.012 & && 0.001 & 0.006 \\
$\delta$&0.2 &  & 0.048 & 0.003 & & 0.044 & 0.002 & && 0.038 & 0.001 \\
$d$&0.7 &  & 0.055& 0.003 & & 0.030& 0.001 &&& 0.026 & 0.0008 \\
$\eta$&1 &  & 0.084 & 0.019 & & 0.057 & 0.019 & && 0.025 & 0.018 \\
 \hline
\end{tabular}
\end{center}
\end{small}
\end{table}

\begin{table}[t] 
\begin{center}
\caption{Empirical rejection rates of the score test for the TV-HGARCH model based on 1000 replications for two significance level 0.05 and 0.10. Also $\eta=0$ corresponds to the size and  $\eta>0$ correspond to the power of the test.}\
\end{center}

\begin{small}
\begin{center}
\begin{tabular}{c c c c c c c c c c c }
\hline \hline
   & &  \multicolumn{2}{c}{n=300} & & \multicolumn{2}{c}{n=500} & & \multicolumn{2}{c}{n=1000} \\
 \cline{3-4} \cline{6-7} \cline{9-10}
 $\eta$   & & 0.05&  0.10 & &  0.05& 0.10 & &  0.05& 0.10 &\\
\hline
$0$ &  & 0.069 & 0.112 & & 0.057 & 0.108 & & 0.049 & 0.095 \\
$0.4$ &  & 0.247 & 0.421 & & 0.547 & 0.731 & & 0.838 & 0.913\\
$1.5$ &  & 0.290 & 0.474 & & 0.578 & 0.759 & & 0.889 & 0.958 \\
$3$ &  & 0.294 & 0.492 & & 0.613 & 0.791 & & 0.915 & 0.966\\
 \hline
\end{tabular}
\end{center}
\end{small}
\end{table}

\section{ VaR Forecasting}
 In order to investigate the ability of the TV-HGARCH model in forecasting the future behavior of the volatilities, we study the VaR  forecasts. The one-day-ahead VaR with  probability $\rho$, $VaR(\rho)$,  is calculated by $VaR_t(\rho)=F^{-1}(\rho)\sigma_t,$  where $F^{-1}(\rho)$ is the inverse distribution of standardized observation $(y_t/\sigma_t)$ and $\sigma_t=\sqrt{V(y_t|\Upsilon_{t-1})}$. Due to the importance of VaR in management risk, the accuracy of the VaR forecasts from different models is evaluated based on  some likelihood ratio (LR) tests (Ardia, 2009; Brooks and Persand, 2000).

\subsubsection*{Unconditional Coverage test}
The Kupiec test (Kupiec, 1995), also known as the unconditional coverage (UC) test, is designed to test whether  VaR forecasts  cover the pre-specified probability.  If the actual loss exceeds the VaR forecasts, this is termed an ``exception,'' which is a Bernoulli random variable with probability $\xi$. The null hypothesis of the UC test is $H_0:\xi=\psi$. Then the LR statistic of the unconditional coverage ($LR_{UC}$) is defined as 
$$
LR_{UC}=-2\log(\dfrac{\rho^n(1-\rho)^{T-n}}{\hat{\xi}^n(1-\hat{\xi})^{T-n}}).$$
Where $T$ is the number of the forecasting samples, $n$ is the number of the exceptions and $\hat{\xi}=\dfrac{n}{T}$ is the MLE of the $\xi$ under $H_1$. Then  under $H_0$ the $LR_{UC}$ is asymptotically distributed as a $\chi^2$ random variable with one degree of freedom.
\subsubsection*{Independent Test}
If the volatilities are low in some periods and high in others, the forecasts should respond to this clustering event. It means that,  the exceptions should be spread over the entire sample period independently and do not appear in clusters (Sarma et al., 2003). Christoffersen (1998) designed an independent (IND) test to check the clustering of the exceptions. The null hypothesis of the IND test assumes that the probability of an exception on a given day t is not influenced by what happened the day before. Formally, $H_0:\xi_{10}=\xi_{00}$, where $\xi_{ij}$ denotes that the probability of an $i$ event on day $t-1$ must be followed by a $j$ event on day $t$ where $i,j=0,1$. The LR statistic of the IND test ($LR_{IND}$) can be obtained as 
$$
LR_{IND}=-2\log(\dfrac{\hat{\xi}^n\hat{(1-\xi)^{T-n}}}{\hat{\xi}_{01}^{n_{01}}(1-\hat{\xi}_{01})^{n_{00}}\hat{\xi_{11}}^{n_{11}}(1-\hat{\xi_{11}})^{n_{10}}}).$$
Where $n_{ij}$ is the number of observations with value $i$ followed by value $j$ ($i,j=0,1$),  $\xi_{01}=\dfrac{n_{01}}{n_{00}+n_{01}}$ and $\xi_{11}=\dfrac{n_{11}}{n_{10}+n_{11}}$. Under $H_0$, the $LR_{UC}$ is asymptotically distributed as a $\chi^2$ random variable with one degree of freedom. 
\subsubsection*{Conditional Coverage test}
Also Christoffersen (1998) proposed a joint test: the conditional coverage (CC) test, which combines the properties of both the UC and IND tests. The null hypothesis of the CC test checks both the exception cluster and consistency of the exceptions with VaR confidence level. The null   hypothesis of the  test is $H_0:\xi_{01}=\xi_{11}=\rho$. The LR  statistic of the CC test ($LR_{CC})$ is obtained as
$$LR_{CC}=-2\log(\dfrac{\rho^n(1-\rho)^{T-n}}{\hat{\xi}_{01}^{n_{01}}(1-\hat{\xi}_{01})^{n_{00}}\hat{\xi_{11}}^{n_{11}}(1-\hat{\xi_{11}})^{n_{10}}}).$$
Under $H_0$, $LR_{CC}$ is asymptotically distributed as a $\chi^2$ random variable with two degrees of freedom. It is a summation of two separate statistics, $LR_{UC}$ and $LR_{IND}$.

\section{Empirical Data}
In this section, we apply the TV-HGARCH(1,d,1), HGARCH(1,d,1)  and FIGARCH(1,d,1) models on the daily percentage log-returns of the \textit{S}\&\textit{P}500 indices from February 17, 2009 to January 30, 2015 (1500 observations). Figure 1 presents the  time plot of data.
Some descriptive statistics of the \textit{S}\&\textit{P}500 indices are listed in Table 3. We observe the negative skewness and excess kurtosis of these returns. To compare the empirical performance of the models from both fitting and forecasting the whole sample is divided into two parts. The first part contains 1,000 observations and is used as in-sample data to conduct fitting and  the second part is used as out-of-sample data to evaluate model forecasting.  The  models are then applied to the first part of data. The MLE values  are reported in Table 4. 
Also the score test in section 4 is performed and the value $\lambda_s=6.44$  is obtained, so the critical value $3.84$  shows that at 5$\%$ significance level the data possesses a time-varying amplitude. To evaluate the performance of the  models in computing true conditional variances that are measured by squared returns, we considered the root mean squared error (RMSE) and the log likelihood value (LLV) for in-sample and out-of-sample data. As out-of-sample performance, the one-day-ahead forecasts are computed using estimated models. The results are given in Table 5.  It is observed  that the TV-HGARCH model has the best performance. The HGARCH model outperforms the FIGARCH model, and has a lower RMSE and a higher LLV. To clarify the out-performance of the TV-HGARCH model, we plot the forecasting conditional variances and true conditional variances for some of the data in Figure 2. It can be seen  that the TV-HGARCH follows the shocks very well. Figure 3 shows the absolute forecasting errors between different models and the true conditional variances for some of the data, it can be observed that the TV-HGARCH model has the  smallest absolute error. 
Based on the out-of-sample data, one-day-ahead VaR forecasts at a level risk of $\rho=0.05,0.10$ for the models are calculated and the accuracy tests are performed. The results are reported in Table 6. The first  and second rows show the number of expected exceptions (Ex.e) and empirical exceptions (Em.e) respectively. It can be seen that the Em.e for the TV-HGARCH model is closer to the Ex.e than HGARCH and FIGARCH models; in this respect also the HGARCH make better results than the FIGARCH. For VaR(0.05) at 5$\%$ significance level,  the TV-HGARCH  model passes all the tests while the HGARCH model passes IND and CC tests and FIGARCH model  passes only IND test. Also for VaR(0.10), all models pass only IND test  but the TV-HGARCH model has the smallest $LR_{UC}$ and $LR_{CC}$. Hence, the results indicate that the TV-HGARCH model produces the most accurate VaR forecasts. Also the HGARCH model outperforms the FIGARCH model.
\begin{table}[t]
\caption{Descriptive statistics of  \textit{S}\&\textit{P}500 daily log-returns}
\begin{small}
\begin{center}
\begin{tabular}[t]{c c c c c c c} \hline
series&Mean&Std.dev&Minimum&Maximum&Skewness&Kurtosis \\ \hline
\textit{S}\&\textit{P}&0.062&1.114&-6.896&6.837&-0.148&4.564\\
 \hline
\end{tabular}
\end{center}
\end{small}
\end{table}

\begin{figure}[!hbtp]
\begin{center}
\includegraphics[width=17cm]{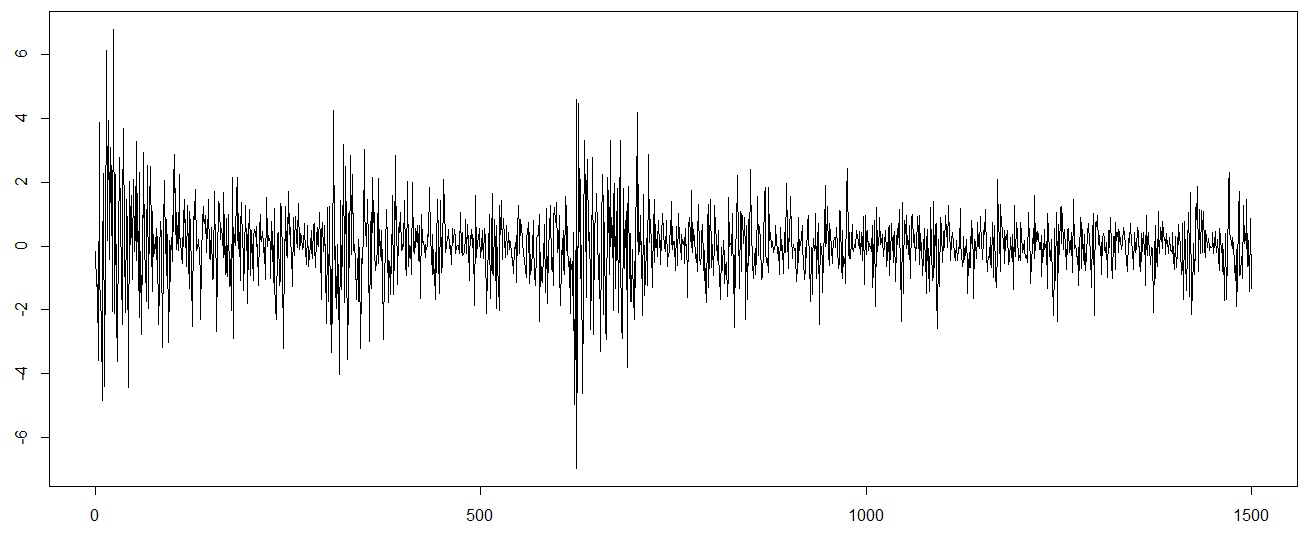}
\caption{Percentage log returns of  \textit{S}\&\textit{P}500 daily log-returns.}
\end{center}
\end{figure}

\begin{center}
\begin{figure}[!btp]
\includegraphics[width=17cm]{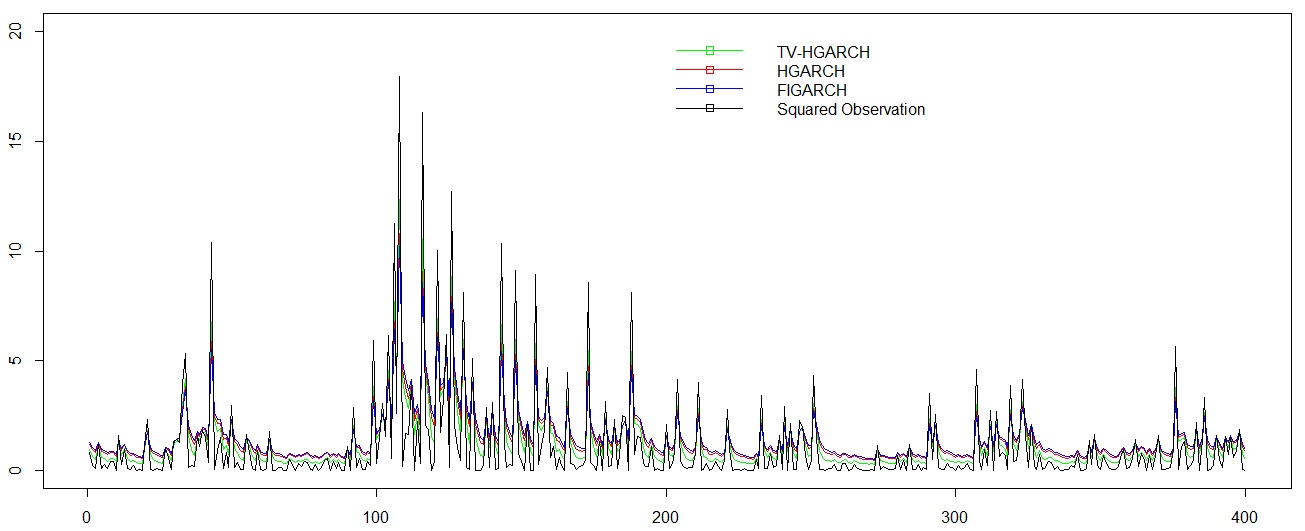}
\caption{Squared returns and forecasting conditional variances with TV-HGARCH, HGARCH  and FIGARCH models for  some of \textit{S}\&\textit{P}500  daily log-returns.} 
\end{figure}
\end{center}

\begin{center}
\begin{figure}[!btp]
\includegraphics[width=17cm]{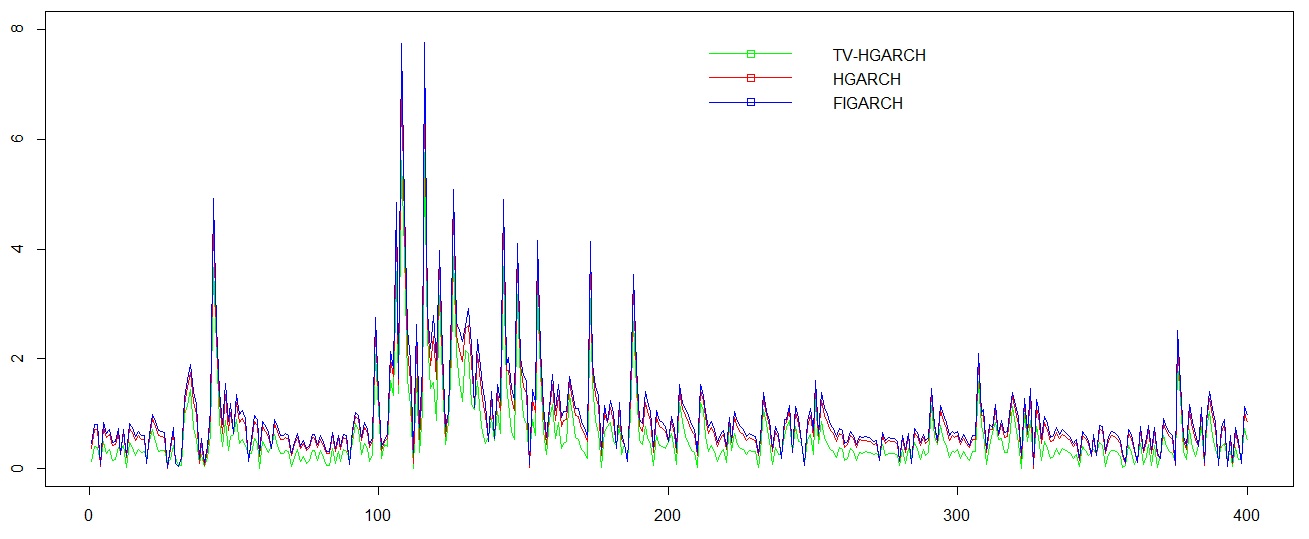}
\caption{Absolute forecasting errors between squared returns and forecasting conditional variances with TV-HGARCH,  HGARCH and FIGARCH models for  some of \textit{S}\&\textit{P}500  daily log-returns.}
\end{figure}
\end{center}
\begin{table}[!hpbt]
\caption{Maximum likelihood estimates of TV-HGARCH, HGARCH and  FIGARCH models on {\textit{S}\&\textit{P}500} daily log-returns. }
\begin{small}
\begin{center}
\begin{tabular}{c c c c c c c c c c c}
\hline
 & TV-HYGARCH &HGARCH & FIGARCH  \\
\hline
$\gamma$&0.179 &0.237&0.237 \\
$\beta$  &0.340&0.455 & 0.505\\
$\delta$ &0.315&0.315&0.315\\
$d$&0.550&0.567&0.505\\
 $\eta$ & 2.444&-&-&\\
 $w$ &-&0.972&-&\\
 \hline
\end{tabular}
\end{center}
\end{small}
\end{table}

\begin{table}[t] 
\caption{Measures of performance  of TV-HGARCH, HGARCH and  FIGARCH models on {\textit{S}\&\textit{P}500} daily log- returns}
\begin{small}
\begin{center}
\begin{tabular}{c c c c c c c c c  }
\hline \hline
   & &  \multicolumn{2}{c}{In-Sample} & & \multicolumn{2}{c}{Out-of-Sample}  \\
 \cline{3-4} \cline{6-7} 
 Model   & & RMSE&  LLV & &  RMSE& LLV \\
\hline
TV-HGARCH &  & 1.406& -1187.4 & & 0.431 & -447.9 \\
HGARCH &  & 1.808& -1294.3 & & 0.637 & -514.6\\
FIGARCH  & & 1.988 & -1325.6& & 0.710 & -530.5\\
 \hline
\end{tabular}
\end{center}
\end{small}
\end{table}
\begin{table}[t] 
\caption{VaR forecasting for TV-HGARCH, HGARCH and FIGARCH models on {\textit{S}\&\textit{P}}500 daily log-returns at level $\rho=0.05,0.10$.}
\begin{small}
\begin{center}
\begin{tabular}{c c c c c c c c c c c }
\hline \hline
   &   \multicolumn{2}{c}{TV-HGARCH} &&  \multicolumn{2}{c}{HGARCH}&&\multicolumn{2}{c}{FIGARCH}   \\
 \cline{2-3} \cline{5-6} \cline{8-9}
  & VaR(0.05)& VaR(0.10) & & VaR(0.05)& VaR(0.10)&& VaR(0.05)& VaR(0.10) \\
\hline
Ex.e &  25 & 50 & & 25  &50& & 25  &50\\
Em.e &  21&33& & 16 & 29&&13&28&  \\
$LR_{UC}$&$0.711^*$ &7.210&& 3.890& 11.371&& 7.298&12.588\\
$LR_{IND}$ &$1.932^*$&$0.954^*$&&$1.125^*$&$0.481^*$&&$0.748^*$&$0.379^*$\\
$LR_{CC}$&$2.642^*$&$8.164$&&$5.013^*$&11.852&&8.047&12.967\\ 
 \hline
\end{tabular}
\end{center}
\end{small}
\begin{small}
Notes: 1. At the 5\% significance level the critical value of the $LR_{UC}$and $LR_{IND}$ is 3.84 and for $LR_{CC}$ is 5.99. 2. * indicates that the model passes the test at 5\% significance level.
\end{small}
\end{table}
{\bf \large  Conclusion:}\\\
HGARCH is a  hyperbolic-memory process. In this study we have generalized  it by introducing TV-HGARCH model to have a better description  of the dynamic  volatilities. Our proposed model exploits a  time-varying amplitude to update the structure of the volatility using logistic function of last observation. We show under some conditions the moments of model are bounded. One of the privilege of this work is implying of score test to check existence of the time-varying structure. Simulation evidences showed that empirical performance of test  is competitive. The empirical example of some periods of \textit{S}\&\textit{P}500 indices showed that the TV-HGARCH model gives better forecasting of volatilities and more accurate VaR  than  HGARCH and FIGARCH.\\

\end{document}